\documentclass[11pt,letterpaper]{amsart}
\usepackage{amsmath,amssymb,amsthm,mathscinet}
\usepackage{enumerate}
\usepackage{multirow,color,colortbl}
\renewcommand{\phi}{\varphi}
\DeclareMathOperator{\Gal}{Gal}
\DeclareMathOperator{\id}{id}
\DeclareMathOperator{\Tr}{Tr}
\DeclareMathOperator{\wt}{wt}
\newcommand{\C}{{\mathbb C}}
\newcommand{\F}{{\mathbb F}}
\newcommand{\Q}{{\mathbb Q}}
\newcommand{\Z}{{\mathbb Z}}
\newcommand{\Fp}{\F_p}
\newcommand{\Fu}{F^*}
\newcommand{\mchars}{\widehat{\Fu}}
\newcommand{\valp}{v_p}
\newcommand{\valtwo}{v_2}
\newcommand{\card}[1]{\left|{#1}\right|}
\newcommand{\conj}[1]{\overline{#1}}
\newcommand{\mins}[1]{\min_{\substack{#1}}}
\newcommand{\floor}[1]{\lfloor{#1}\rfloor}
\newcommand{\ceil}[1]{\lceil{#1}\rceil}
\newtheorem{theorem}{Theorem}[section]
\newtheorem{lemma}[theorem]{Lemma}
\newtheorem{corollary}[theorem]{Corollary}
\newtheorem{conjecture}[theorem]{Conjecture}
\newtheorem{openproblem}[theorem]{Open Problem}
\theoremstyle{definition}
\newtheorem{definition}[theorem]{Definition}
\definecolor{tablecolor}{rgb}{0.85,0.85,0.85}
\title[Weil sums of binomials]{Weil sums of binomials: properties, applications, and open problems}
\author{Daniel J. Katz}
\address{Department of Mathematics, California State University, Northridge, \: United States}\thanks{This paper is based upon work supported in part by the National Science Foundation under Grants DMS-1500856 and CCF-1815487.}
\date{19 November 2018}
\begin{document}
\begin{abstract}
We present a survey on Weil sums in which an additive character of a finite field $F$ is applied to a binomial whose individual terms (monomials) become permutations of $F$ when regarded as functions.
Then we indicate how these Weil sums are used in applications, especially how they characterize the nonlinearity of power permutations and the correlation of linear recursive sequences over finite fields.
In these applications, one is interested in the spectrum of Weil sum values that are obtained as the coefficients in the binomial are varied.
We review the basic properties of such spectra, and then give a survey of current topics of research: Archimedean and non-Archimedean bounds on the sums, the number of values in the spectrum, and the presence or absence of zero in the spectrum.
We indicate some important open problems and discuss progress that has been made on them.
\end{abstract}
\maketitle

\section{Weil sums of binomials}\label{Wally}

A {\it Weil sum} is a finite field character sum where the character is evaluated with a polynomial argument.
We are interested in additive character Weil sums, that is, sums of the form
\[
\sum_{x \in F} \phi(f(x)),
\]
where $F$ is a finite field, $\phi\colon F \to \C$ is an additive character, and $f(x)$ is a polynomial in $F[x]$.
We typically use the {\it canonical additive character}, $\psi \colon F \to \C$ with $\psi(x)=e^{2\pi i \Tr(x)/p}$, where $\Tr \colon F \to \Fp$ is the absolute trace, because for any additive character $\phi \colon F \to \C$ there is some $c \in F$ such that $\phi(x)=\psi(c x)$ for all $x \in F$, so that the above sum can be recast into the form
\[
\sum_{x \in F} \psi(g(x)),
\]
where $g(x)=c f(x) \in F[x]$.
In this paper, we are specifically interested in Weil sums with binomial arguments, that is, sums of the form
\begin{equation}\label{Edgar}
\sum_{x \in F} \psi(a x^d + b x^e),
\end{equation}
where $a,b \in \Fu$ and $d\not=e$.
Such sums have been studied in \cite{Mordell,Davenport-Heilbronn,Karatsuba,Carlitz-78,Carlitz-79,Coulter,Cochrane-Pinner-03,Cochrane-Pinner-11}.
For example, we see that the Kloosterman sum,
\[
\sum_{x \in \Fu} \psi(a x+ b x^{-1})
\]
is one less than what one obtains in \eqref{Edgar} in the special case where $d=1$ and $e=\card{F}-2$.

We are particularly interested in the case where the exponents that appear in our binomial are invertible modulo $\card{\Fu}$.
\begin{definition}[Invertible Exponent]
If $F$ is a finite field, we say that a positive integer $d$ is an {\it invertible exponent over $F$} (or just say that $d$ is {\it invertible over $F$}) to mean that $\gcd(d,\card{\Fu})=1$.
\end{definition}
The map $x \mapsto x^d$ is a permutation of the finite field $F$ if and only if $d$ is invertible over $F$, and we have a special name for such permutations.
\begin{definition}[Power Permutation]\label{Jane}
If $F$ is a finite field, a {\it power permutation} is a mapping from $F$ to itself of the form $x \mapsto x^d$ for some $d$ that is invertible over $F$.
\end{definition}
If $d$ and $e$ in \eqref{Edgar} are invertible over $F$, we may use this fact to reparameterize the Weil sum by replacing $x$ with $a^{-1/d} x^{1/e}$ (where inversion in exponents is modulo $\card{\Fu}$) to obtain
\[
\sum_{x \in F} \psi(x^{d/e}+b a^{-e/d} x),
\]
so it suffices to study sums in the following standard form.
\begin{definition}[$W_{F,d}$]\label{William}
If $F$ is a finite field with canonical additive character $\psi$, and $d$ is invertible over $F$, then $W_{F,d}$ always denotes the function from $F$ to $\C$ with
\[
W_{F,d}(a)= \sum_{x \in F} \psi(x^d - a x),
\]
for each $a \in F$.
\end{definition}
For fixed $F$ and $d$, we are interested in the values of $W_{F,d}(a)$ as $a$ runs through $\Fu$.
\begin{definition}[Weil Spectrum]
If $F$ is a finite field and $d$ is an invertible exponent over $F$, then the {\it Weil spectrum} associated with $F$ and $d$ is the set $\{W_{F,d}(a): a \in \Fu\}$.
\end{definition}
Readers familiar with the Walsh spectrum of a power permutation (introduced later in this paper in Section \ref{Persephone}) should note that the Walsh spectrum includes the value of $W_{F,d}(0)$ while the Weil spectrum defined here does not.
When one does consider $W_{F,d}(0)$, the polynomial in the definition of $W_{F,d}$ is a monomial, and the fact that $x \mapsto x^d$ is a permutation of $F$ makes $W_{F,d}(0)=0$.

In this survey we first review the applications of the Weil sums $W_{F,d}(a)$ to nonlinearity of finite field permutations (Section \ref{Persephone}), correlation of pseudorandom sequences (Section \ref{Martha}), and coding theory and finite geometry (Section \ref{Elijah}).
We then review their basic properties using the Galois theory of finite fields (Section \ref{Bartholomew}), the Galois theory of cyclotomic extensions of $\Q$ (Section \ref{Cecilia}), and power moments (Section \ref{Daria}).
Finally we survey the current state of knowledge concerning these sums, including bounds on them, both Archimedean (Section \ref{Edwin}) and non-Archimedean (Section \ref{Felix}), the number of distinct values in the Weil spectrum (Section \ref{Gerald}), and the presence or absence of zero in the Weil spectrum (Section \ref{Nigel}).
These last four sections state some of the principal open problems in the field, and indicate progress toward their resolution.

\section{Nonlinearity of power permutations}\label{Persephone}

Power permutations, introduced in Definition \ref{Jane}, are of interest in cryptography, for example, in the S-boxes of symmetric key cryptosystems.
We are interested in how far our power permutation is from being linear.
To this end, we define a transform that measures nonlinearity.
\begin{definition}[Walsh Transform]\label{Eric}
Let $F$ be a finite field of characteristic $p$ and $f\colon F \to F$.
Then the {\it Walsh transform} of $f$ is the function $W_f \colon F^2 \to \C$ given by
\[
W_f(b,a)=\sum_{x \in F} \psi(b f(x)-a x),
\]
where $\psi(y)=e^{2\pi i \Tr(y)/p}$ is the canonical additive character of $F$, with $\Tr$ the absolute trace of $F$.
\end{definition}
When $p=2$ we have $\psi(y)=(-1)^{\Tr(y)}$ so $W_f(b,a)$ measures the number of agreements minus the number of disagreements between the functions $x \mapsto \Tr(b f(x))$ and the $\F_2$-linear functions $x\mapsto \Tr(a x)$.
In other characteristics $p$, the meaning of the Walsh transform values is not as obvious, but a large magnitude value still indicates a high degree of resemblance.
Since we are trying to avoid any resemblance, we should examine all the values $W_f(b,a)$ except those with $b=0$, for $W_f(0,a)$ tells us nothing about $f$: it is trivially $\card{F}$ when $a=0$ and trivially $0$ otherwise.
\begin{definition}[Walsh Spectrum]
Let $F$ be a finite field and $f\colon F \to F$ with Walsh transform $W_f$.
The {\it Walsh spectrum of f} is the set of values $\{W_f(b, a): b \in \Fu , a \in F \}$.
\end{definition}
A highly nonlinear function $f$ is one where all the values in the Walsh spectrum have small magnitude.
When $d$ is an invertible exponent over our field $F$, and $f$ is the power permutation $f(x)=x^d$, then we have Walsh transform
\begin{equation}\label{Francis}
W_f(b,a)=\sum_{x \in F} \psi(b x^d-a x).
\end{equation}
When $b\not=0$, the reparameterization described in Section \ref{Wally} shows that $W_f(b,a)=W_{F,d}(b^{-1/d} a)$.
Thus we see that the Walsh spectrum for our power permutation is given by the values of $W_{F,d}(a)$ as $a$ runs through $F$, that is, the Walsh spectrum is just the Weil spectrum along with $W_{F,d}(0)=0$.
So a highly nonlinear power permutation $x\mapsto x^d$ of the finite field $F$ is one where $|W_{F,d}(a)|$ is small for all $a \in F$.
We shall see later (in Lemma \ref{Jacob}) that it is not possible to have all $|W_{F,d}(a)|$ less than or equal to $\sqrt{\card{F}}$, but Table \ref{Nancy} in Section \ref{Gerald} exhibits cases in which we can make every $|W_{F,d}(a)|$ value less than or equal to $\sqrt{p \card{F}}$, where $p$ is the characteristic of $F$.

\section{Crosscorrelation of m-sequences and Gold sequences}\label{Martha}

Linear recursive sequences over finite fields have been employed to produce pseudorandom sequences as early as the 1950's by Gilbert, Golomb, Welch, and Zierler (see the Preface of \cite{Golomb} and the Historical Introduction of \cite{Golomb-Gong} for brief historical overviews).
Much of the mathematical theory of these sequences had been worked out somewhat earlier: see \cite{Carmichael,Engstrom,Ward-31,Ward-33,Hall}.
From this theory, one can show that a sequence of elements from the field $\Fp$ of prime order $p$ satisfying a linear recursion whose minimal polynomial is of degree $n$ over $\Fp$ has a period less than or equal to $p^n-1$, with equality if and only if the minimal polynomial is a primitive irreducible.
In this case, the sequence is called a {\it maximum length linear feedback shift register sequence}, or a {\it maximal linear sequence}, or an {\it m-sequence}.

Golomb had given criteria for evaluating the pseudorandomness of binary sequences \cite[p.~25]{Golomb}, which can be generalized to nonbinary sequences \cite[\S5.1.2]{Golomb-Gong}.
One of these criteria is that the sequence should have low autocorrelation at nontrivial shifts; this property makes these sequences useful for applications to remote sensing and communications.
When communications systems have many users, it also becomes useful to have families of such sequences that have the additional property of low crosscorrelation between any pair of sequences in the family.
The m-sequences have very low autocorrelation, and pairs of m-sequences with the same period can be found that have relatively low crosscorrelation, and these in turn can be used to construct larger families of sequences (Gold sequences) in which every pair of sequences has low periodic crosscorrelation.

A sequence with terms from the prime field $\Fp$ is called a {\it $p$-ary sequence}.
It is natural to index the terms of any $p$-ary m-sequence of period $p^n-1$ using $\Z/(p^n-1)\Z$.
If we fix $\alpha$, a primitive element of order $p^n-1$ in the finite field $F$ of order $p^n$, and let $\Tr \colon F \to \Fp$ be the absolute trace, then any $p$-ary m-sequence of period $p^n-1$ can be written as
\begin{equation}\label{Albert}
(\Tr(\alpha^{d j+s}))_{j \in \Z/(p^n-1)\Z}
\end{equation}
for some integers $s$ and $d$ with $\gcd(d,p^n-1)=1$.
Let us take the m-sequence
\begin{equation}\label{Rachel}
(\Tr(\alpha^j))_{j \in \Z/(p^n-1)\Z}
\end{equation}
with $d=1$ and $s=0$ as our reference sequence.
Then for any $d$ with $\gcd(d,p^n-1)=1$, the m-sequence
\begin{equation}\label{Desdemona}
(\Tr(\alpha^{d j}))_{j \in \Z/(p^n-1)\Z}
\end{equation}
is obtained from our reference sequence by taking every $d$th element (proceeding cyclically according to our indexing): we call this new sequence the {\it decimation by $d$} of the reference sequence.
And then we see that a generic m-sequence \eqref{Albert} of period $p^n-1$ is obtained from our reference sequence by decimating by $d$ and then shifting (cyclically) by $s$ places.
When $d$ is a power of $p$ modulo $p^n-1$, the decimated sequence is identical to the reference sequence, and we say that such a decimation is {\it trivial}.

We shall now study the correlation properties of m-sequences, after setting down the basic definitions.
\begin{definition}[Periodic Crosscorrelation]
Let $\ell$ be a positive integer, and let $f=(f_j)_{j \in \Z/\ell\Z}$ and $g=(g_j)_{j \in \Z/\ell\Z}$ be $p$-ary sequences whose periods are divisors of $\ell$.
For any $s \in \Z/\ell\Z$, the {\it periodic crosscorrelation of $f$ with $g$ at shift $s$}, denoted $C_{f,g}(s)$, is defined to be
\[
C_{f,g}(s)=\sum_{j \in \Z/\ell\Z} e^{2\pi i (f_{j+s}-g_j)/p}.
\]
\end{definition}
We are interested in the crosscorrelation values at all the various shifts.
\begin{definition}[Crosscorrelation Spectrum]
Let $\ell$ be a positive integer, and let $f=(f_j)_{j \in \Z/\ell\Z}$ and $g=(g_j)_{j \in \Z/\ell\Z}$ be $p$-ary sequences whose periods are divisors of $\ell$.
The {\it crosscorrelation spectrum} of the sequence pair $(f,g)$ is $\{C_{f,g}(s): s \in \Z/\ell\Z\}$.
\end{definition}
It turns out that crosscorrelation values are character sums when $f$ and $g$ are m-sequences.
Let $g$ be our reference m-sequence \eqref{Rachel} and $f$ be \eqref{Desdemona} obtained by decimating $g$ by $d$.
We let $\psi\colon F \to \C$ be the canonical additive character, $\psi(x)=e^{2\pi i \Tr(x)/p}$, of our field $F$ of characteristic $p$ and order $p^n$.
Then
\begin{align}
C_{f,g}(s)
& =\sum_{j \in \Z/(p^n-1)\Z} \psi(\alpha^{d (j+s)}-\alpha^j) \\
& =\sum_{j \in \Z/(p^n-1)\Z} \psi(\alpha^{d j}-\alpha^{j-s}) \\
& =\sum_{x \in \Fu}\psi(x^d-\alpha^{-s} x) \\
& =-1+W_{F,d}(\alpha^{-s}),
\end{align}
where $W_{F,d}$ is the Weil sum from Definition \ref{William}.
Thus the crosscorrelation spectrum for an m-sequence of period $p^n-1$ and its decimation by $d$ is given by the values of $-1+W_{F,d}(a)$ as $a$ runs through $\Fu$, that is, the crosscorrelation spectrum is obtained from the Weil spectrum by subtracting $1$ from every value.
Note that when the decimation $d=1$, we are crosscorrelating our m-sequence with itself (this is called {\it autocorrelation}).
In this case the autocorrelation values $-1+W_{F,1}(a)$ all equal $-1$, except when $a=1$ (which corresponds to zero shift, and has autocorrelation $\card{F}-1=p^n-1$).
The property that autocorrelation at every nontrivial shift has a value of $-1$ is called the {\it ideal two-level autocorrelation} property of m-sequences, and is useful in applications involving synchronization.

If we continue to identify $g$ and $f$ with the m-sequences from \eqref{Rachel} and \eqref{Desdemona}, then for each $b \in F$, we can define the sequence
\[
h^{(b)}= (\Tr(\alpha^{d j}-b\alpha^j))_{j \in \Z/(p^n-1)\Z},
\]
and we introduce the symbol $\infty$, and define
\[
h^{(\infty)} = g = (\Tr(\alpha^j))_{j \in \Z/(p^n-1)\Z}.
\]
This gives a family $\{h^{(b)}: b \in F\cup\{\infty\}\}$ of $p^n+1$ sequences of period $p^n-1$ known as a family of {\it Gold sequences} \cite{Gold-67}.\footnote{Some authors only refer to the family so produced as a family of Gold sequences if the decimation $d$ is of a particular form originally selected by Gold, which produces favorable correlation spectra; but in the final portion of Gold's paper \cite{Gold-67}, the construction is presented in a way that a reader could straightforwardly generalize to an arbitrary decimation.}
For any $b \in F\cup\{\infty\}$, we can compute the autocorrelation function of the Gold sequence $h^{(b)}$ to obtain
\[
C_{h^{(b)},h^{(b)}}(s) = \begin{cases}
p^n-1 & \text{if $s=0$,} \\
-1 & \text{if $b=\infty$ and $s\not=0$,} \\
-1+W_{F,d}\left(\frac{b(\alpha^s-1)}{(\alpha^{d s}-1)^{1/d}}\right) & \text{if $b\not=\infty$ and $s\not=0$,}
\end{cases}
\]
where $1/d$ indicates the multiplicative inverse of $d$ modulo $\card{\Fu}$.
So the autocorrelation values at nonzero shifts that one obtains among the various Gold sequences in our family are given by values of $-1+W_{F,d}(a)$ as $a$ runs through $F$ (recall that $W_{F,d}(0)=0$ while considering the penultimate case).
Now we suppose that $b$ and $c$ are two distinct elements of $F \cup\{\infty\}$, and compute the crosscorrelation between $h^{(b)}$ and $h^{(c)}$ at shift $s \in \Z/(p^n-1)\Z$ to be
\[
C_{h^{(b)},h^{(c)}}(s) = \begin{cases}
-1+W_{F,d}(b+\alpha^{-s}) & \text{if $b\not=\infty$ and $c=\infty$,} \\
-1+W_{F,d}(c+\alpha^s) & \text{if $b=\infty$ and $c\not=\infty$,} \\
-1 & \text{if $b,c\not=\infty$, and $s=0$,} \\
-1+W_{F,d}\left(\frac{b\alpha^s-c}{(\alpha^{d s}-1)^{1/d}}\right) & \text{if $b,c\not=\infty$ and $s\not=0$,}
\end{cases}
\]
so we see that the crosscorrelation values one obtains among pairs of distinct Gold sequences in our family are given by values of $-1+W_{F,d}(a)$ as $a$ runs through $F$ (and again, recall that $W_{F,d}(0)=0$ while considering the penultimate case).

In applications involving Gold sequences, we want the magnitudes of autocorrelations at nonzero shifts and of crosscorrelations at all shifts to be small.
This is equivalent to saying that we want $|-1+W_{F,d}(a)|$ to be small for all $a \in F$, which is essentially the same goal we stated at the end of Section \ref{Persephone} for obtaining highly nonlinear power permutations (which was that $|W_{F,d}(a)|$ be small for all $a \in F$).

\section{Other applications}\label{Elijah}

Suppose that $F$ is the finite field of characteristic $p$ and order $q=p^n$ with primitive element $\alpha$, and let $d$ be a positive integer with $\gcd(d,q-1)=1$.
If $d \equiv 1 \pmod{p-1}$ and $d$ is not congruent to a power of $p$ modulo $q-1$, then the values of the Weil sum $W_{F,d}(a)$ in Definition \ref{William} as $a$ runs through $F$ determine the weight distribution of the $p$-ary cyclic error-correcting code ${\mathcal C}_{1,d}$ of length $p^n-1$ and dimension $2 n$ whose check polynomial is the product of the minimal polynomials of $\alpha^{-1}$ and $\alpha^{-d}$ over $\Fp$ (see \cite[\S A.3]{Katz-12}).
This code is the sum of two simplex codes, ${\mathcal C}_1$ and ${\mathcal C}_d$, whose respective check polynomials are the minimal polynomials of $\alpha^{-1}$ and $\alpha^{-d}$ over $\Fp$.
The nonzero words of these simplex codes all have weights of $(p-1) q/p$, while the words of ${\mathcal C}_{1,d}\smallsetminus ({\mathcal C}_1 \cup {\mathcal C}_d)$ have weights of the form
\begin{equation}\label{Ulrich}
\frac{(p-1) (q-W_{F,d}(a))}{p},
\end{equation}
where $a$ runs through $\Fu$.
Since $W_{F,d}(0)=0$, we can summarize the above facts by saying that the nonzero words of ${\mathcal C}_{1,d}$ have weights given by \eqref{Ulrich} as $a$ runs through $F$.

When $p=2$, the values of the Weil sum $W_{F,d}(a)$ in Definition \ref{William} also have an application in finite geometry (see \cite[\S A.4]{Katz-12}).
They give the cardinalities of the intersections of hyperplanes with certain constructions described by Games in \cite{Games-m-sequences,Games-quadrics}.

\section{Galois Action of the Finite Field}\label{Bartholomew}

In this section we shall give known results about the Weil sum $W_{F,d}(a)$ from Definition \ref{William} arising from the Galois theory of the finite field $F$.
\begin{lemma}
If $\sigma$ is any automorphism of the finite field $F$ and $a \in F$, then $W_{F,d}(\sigma(a))=W_{F,d}(a)$.  In other words, if $F$ is of characteristic $p$, then $W_{F,d}(a^{p^k})=W_{F,d}(a)$ for every $a \in F$ and every integer $k$.
\end{lemma}
\begin{lemma}\label{Fiona}
If $F$ is of characteristic $p$ and order $q$, and $e \equiv p^k d \pmod{q-1}$ for some integer $k$, then $W_{F,e}(a)=W_{F,d}(a)$ for every $a \in F$.
\end{lemma}
See \cite[Lemmata 3.1--3.2]{Aubry-Katz-Langevin-15} for proofs of these two lemmata: the equivalent results in terms of crosscorrelation of m-sequences are well-known and go back at least as far as \cite[Theorem 3.1(d)--(e)]{Helleseth-76}.
Lemma \ref{Fiona} shows that if $d$ is a power of $p$ modulo $q-1$, then $W_{F,d}(a)=W_{F,1}(a)=\sum_{x \in F} \psi((1-a) x)$, so that we have a Weil sum of a monomial, which prompts the following definition.
\begin{definition}
If $F$ is a finite field of characteristic $p$ and order $q$, a positive integer $d$ congruent to a power of $p$ modulo $q-1$ is said to be a {\it degenerate exponent over $F$}.
\end{definition}
It is straightforward to determine which fields support nondegenerate invertible exponents.
\begin{lemma}\label{Hubert}
Let $F$ be a finite field.
If $\card{F} \leq 4$, every invertible exponent over $F$ is degenerate.
If $\card{F} > 4$, then there is a nondegenerate invertible exponent over $F$, for example, $\card{F}-2$.
\end{lemma}
The values of a degenerate Weil sum are very straightforward to calculate.
\begin{corollary}\label{Deidre}
If $d$ is degenerate over $F$, then
\[
W_{F,d}(a) = \begin{cases}
\card{F} & \text{if $a=1$,} \\
0 & \text{otherwise.}
\end{cases}
\]
\end{corollary}
The equivalent result in terms of crosscorrelation is \cite[Theorem 3.1(g)]{Helleseth-76}.
Using reparameterizations akin to those discussed in Section \ref{Wally}, we have the following result.
\begin{lemma}
If $d$ is invertible over a finite field $F$, and if $1/d$ denotes the inverse of $d$ modulo $\card{\Fu}$, then $W_{F,1/d}(a)=W_{F,d}(a^{1/d})$ for every $a \in F$.
\end{lemma}
This means that $\{W_{F,1/d}(a): a \in S\}=\{W_{F,d}(a): a \in S\}$ for both $S=F$ and $S=\Fu$.
Lemma \ref{Fiona} also gives cases where two exponents produce equivalent Weil sums, so we formulate the following definition.
\begin{definition}\label{Elizabeth}
If $d$ and $e$ are invertible exponents over a finite field $F$ of characteristic $p$ such that either $d \equiv p^k e \pmod{\card{\Fu}}$ or $d \equiv p^k e^{-1} \pmod{\card{\Fu}}$ for some integer $k$, then $d$ and $e$ are said to be {\it equivalent exponents over $F$}.
\end{definition}
The following result is the compositum of Lemma 3.3 and Corollary 3.4 of \cite{Aubry-Katz-Langevin-15}; note that we make some sign changes here to reflect the different sign convention we use in the definition of our Weil sum of a binomial.
\begin{lemma}
Let $F$ be a finite field of characteristic $p$, and let $E$ be an extension of $F$ with $[E:F]$ a power of a prime $\ell$ different from $p$.  If $1/d$ denotes the multiplicative inverse of $d$ modulo $p-1$, then
\[
W_{E,d}(a) \equiv W_{F,d}([E:F]^{1-1/d} a) \pmod{\ell}
\]
for every $a \in F$.  In particular, if $d$ is degenerate over $F$, then
\[
W_{E,d}(a) \equiv \begin{cases}
\card{F} \pmod{\ell} & \text{if $a=1$} \\
0 \pmod{\ell} & \text{if $a \in F\smallsetminus\{1\}$.}
\end{cases}
\]
\end{lemma}

\section{Galois Action of the Cyclotomic Field}\label{Cecilia}

The canonical additive character of a finite field $F$ of characteristic $p$ maps $F$ into the group of $p$th roots of unity in $\C$.
We set $\zeta_p = e^{2\pi i/p}$, so that our Weil sum $W_{F,d}(a)$ from Definition \ref{William} lies in $\Z[\zeta_p]$.
This is the ring of algebraic integers in the cyclotomic field $\Q(\zeta_p)$.
Thus we can take advantage of the Galois theory of this cyclotomic field to obtain results about our Weil sum.
First we note that an automorphism of the cyclotomic field permutes the spectrum of Weil sum values.
\begin{lemma}
Let $F$ be a finite field of characteristic $p$, and let $d$ be an invertible exponent over $F$.  If $\sigma \in \Gal(\Q(\zeta_p)/\Q)$ with $\sigma(\zeta_p)=\zeta_p^j$, then $\sigma(W_{F,d}(a))=W_{F,d}(j^{1-1/d} a)$ for every $a \in F$, where $1/d$ denotes the multiplicative inverse of $d$ modulo $p-1$.
\end{lemma}
The above lemma is proved in \cite[Theorem 2.1(b)]{Katz-12}, and it has the following corollaries (see \cite[Theorem 3.1(a)]{Helleseth-76} and \cite[Corollary 5.2]{Aubry-Katz-Langevin-15}, respectively).
\begin{corollary}
If $F$ is a finite field and $d$ is an invertible exponent over $F$, then $W_{F,d}(a)$ is a real number for every $a \in F$.
\end{corollary}
\begin{corollary}
Let $F$ be a finite field of characteristic $p$, and let $d$ be an invertible exponent over $F$.  If $A, B \in \Q(\zeta_p)$ are Galois conjugates over $\Q$, then the number of $a \in \Fu$ such that $W_{F,d}(a)=A$ is equal to the number of $a \in \Fu$ such that $W_{F,d}(a)=B$.
\end{corollary}
For certain $d$, the Weil sum values lie in a proper subfield of $\Q(\zeta_p)$.
\begin{lemma}
Let $F$ be a finite field of characteristic $p$, and let $d$ be an invertible exponent over $F$.
Let $K$ be the field generated by adjoining the Weil sum values $\{W_{F,d}(a): a \in F\}$ to $\Q$.
Then $[\Q(\zeta_p):K]$ is the largest divisor $m$ of $p-1$ such that $d \equiv 1 \pmod{m}$; since $\Q(\zeta_p)$ is a cyclic extension of $\Q$, this uniquely identifies $K$.
\end{lemma}
This was proved in \cite[Lemma 5.3]{Aubry-Katz-Langevin-15}, and has the following well known consequence.
\begin{corollary}
Let $F$ be a finite field of characteristic $p$, and let $d$ be an invertible exponent over $F$.
Then $\{W_{F,d}(a): a \in \Fu\} \subseteq \Z$ if and only if $d \equiv 1 \pmod{p-1}$.
\end{corollary}
The version of this for crosscorrelation was proved as \cite[Theorem 4.2]{Helleseth-76}.

\section{Moments}\label{Daria}
Often it is difficult to understand the behavior of the individual values $W_{F,d}(a)$, but it is easier to compute and understand quantities determined by their collective behavior.
\begin{definition}[Power Moment]
If $k$ is a nonnegative integer, then the {\it $k$th power moment of $W_{F,d}$}, denoted $P_{F,d}^{(k)}$, is the sum of the $k$th powers of the Weil spectrum, that is, $P_{F,d}^{(k)}=\sum_{a \in \Fu} W_{F,d}(a)^k$.
\end{definition}
A more general quantity was explored in the context of crosscorrelation in \cite[Theorem 3.2]{Helleseth-76}, and we state an equivalent result in terms of Weil sums here.
\begin{lemma}\label{Ursula}
Let $F$ be a finite field and $d$ an invertible exponent over $F$.
If $k$ is a positive integer and $b_1,\ldots,b_k \in \Fu$, then
\[
\sum_{a \in \Fu} \prod_{j=1}^k W_{F,d}(b_j a) = \frac{\card{F}^2 N_{b_1,\ldots,b_k}-\card{F}^k}{\card{F}-1},
\]
where $N_{b_1,\ldots,b_k}$ is the number of solutions $(x_1,\ldots,x_k) \in F^k$ to the system of equations
\begin{align*}
b_1 x_1 + \cdots + b_k x_k & = 0 \\
x_1^d+\cdots +x_k^d & =0.
\end{align*}
\end{lemma}
Note that $N_{b_1,\ldots,b_k}$ in the preceding lemma is the size of the intersection of a hyperplane and a Fermat variety in $k$-dimensional affine space.
When all the $b_j$'s are equal, we obtain the power moments.
\begin{corollary}
Let $F$ be a finite field and $d$ an invertible exponent over $F$.
If $k$ is a positive integer, then the $k$th power moment of $W_{F,d}$ is
\[
P_{F,d}^{(k)} = \sum_{a \in \Fu} W_{F,d}(a)^k = \frac{\card{F}^2 N^{(k)} - \card{F}^k}{\card{F}-1},
\]
where $N^{(k)}$ is the number of solutions $(x_1,\ldots,x_k) \in F^k$ to the system of equations
\begin{align*}
x_1 + \cdots +x_k & = 0 \\
x_1^d+\cdots +x_k^d & =0.
\end{align*}
\end{corollary}
The first few power moments are of particular interest; see for example \cite[Corollary 2.6]{Katz-15}.
\begin{corollary}\label{Paul}
Let $F$ be a finite field and $d$ an invertible exponent over $F$.
\begin{enumerate}[(i)]
\item $P_{F,d}^{(0)} = \sum_{a \in \Fu} W_{F,d}(a)^0 = \card{F}-1$,
\item\label{Frederick} $P_{F,d}^{(1)} = \sum_{a \in \Fu} W_{F,d}(a)^1 = \card{F}$,
\item\label{Stephen} $P_{F,d}^{(2)} = \sum_{a \in \Fu} W_{F,d}(a)^2 = \card{F}^2$,
\item $P_{F,d}^{(3)} = \sum_{a \in \Fu} W_{F,d}(a)^3 = \card{F}^2 M_1$, and 
\item $P_{F,d}^{(4)} = \sum_{a \in \Fu} W_{F,d}(a)^4 = \card{F}^2 \sum_{a \in \Fu} M_a^2$,
\end{enumerate}
where
\[
M_a = \card{\{x \in F: x^d+(1-x)^d=a\}}.
\]
\end{corollary}
The count $M_a$ in the last corollary is equal to the number of $x \in F$ such that
\[
(x+1)^d-x^d=a,
\]
which one sees by reparameterizing with $-x$ in place of $x$ (and noting that $\gcd(d,\card{F}-1)=1$ makes $d$ odd when $F$ is of odd characteristic).
Thus the power moments of $W_{F,d}$ are intimately connected to the differential spectrum of the power permutation $x\mapsto x^d$ of $F$.
Analysis of the differential spectrum has proved a powerful tool in studying $W_{F,d}$.
For example, Dobbertin's work in \cite{Dobbertin-99-Welch,Dobbertin-99-Niho} was a crucial step in verifying that certain $F$ and $d$ produce very favorable Walsh spectra, which are recorded later in this paper as the seventh and ninth entries of Table \ref{Nancy} in Section \ref{Gerald}.

\section{Archimedean Bounds}\label{Edwin}

In the Introduction we stated that in applications it is desirable to have an invertible exponent $d$ over a finite field $F$ such that $|W_{F,d}(a)|$ is small for all values of $a \in \Fu$.
(Since $W_{F,d}(0)=0$, this is the same as saying that $W_{F,d}(a)$ is small for all $a \in F$.)
Since $W_{F,d}(a)$ is a sum of $\card{F}$ roots of unity in $\C$, it is clear that we always have $|W_{F,d}(a)|\leq \card{F}$.
In fact, the only way to achieve equality is if $d$ is degenerate (see \cite[Theorem 2.1(f)]{Katz-12} and Corollary \ref{Deidre} above).
\begin{lemma}\label{Martin}
Let $F$ be a finite field and $d$ an invertible exponent over $F$.
Then $|W_{F,d}(a)|=\card{F}$ if and only if $d$ is degenerate and $a=1$; otherwise $|W_{F,d}(a)| < \card{F}$.
\end{lemma}
One interesting consequence of this fact, in combination with the values of the first and second power moments in Corollary \ref{Paul}\eqref{Frederick}--\eqref{Stephen}, is that a nondegenerate $W_{F,d}$ must assume at least one positive and at least one negative value \cite[Corollary 2.3]{Aubry-Katz-Langevin-14}.
\begin{lemma}\label{Julie}
If $F$ is a finite field and $d$ is a nondegenerate invertible exponent over $F$, then $W_{F,d}(a) > 0$ for some $a \in F$ and $W_{F,d}(b) < 0$ for some $b$ in $F$.
\end{lemma}
Corollary \ref{Paul}\eqref{Stephen} also implies a useful standard about how small we can hope to make all the Weil sum values.
\begin{lemma}\label{Jacob}
For any finite field $F$ and invertible exponent $d$ over $F$, there is some $a \in F$ such that $|W_{F,d}(a)| > \sqrt{\card{F}}$.
\end{lemma}
This result, combined with the fact that there are known $F$ and $d$ such that $|W_{F,d}(a)| \leq \sqrt{2 \card{F}}$ for all $a \in F$ (for examples, see Table \ref{Nancy} in Section \ref{Gerald}), means that ``small'' should understood as ``not much larger than $\sqrt{\card{F}}$.''

The Weil-Carlitz-Uchiyama bound on character sums \cite{Weil,Carlitz-Uchiyama} can be applied to our Weil sum.
\begin{theorem}
Let $F$ be a finite field and $d$ a nondegenerate invertible exponent over $F$.
Then $|W_{F,d}(a)|\leq (d-1)\sqrt{\card{F}}$ for every $a \in F$.
\end{theorem}
Of course, when $d \geq 1+\sqrt{\card{F}}$, this bound is trivial.
When $d=\card{F}-2$, we see that 
\[
W_{F,\card{F}-2}(a)=1+\sum_{x \in \Fu} \psi(x^{-1}-a x),
\]
that is, one plus a Kloosterman sum.  The results of Weil-Carlitz-Uchiyama \cite{Weil,Carlitz-Uchiyama} can also be applied in this case.
\begin{theorem}
If $F$ is a finite field, then $|W_{F,\card{F}-2}(a)|\leq 2\sqrt{\card{F}}$ for every $a \in F$.
\end{theorem}
Interestingly, if $F$ is of characteristic $2$, then Lachaud and Wolfmann \cite{Lachaud-Wolfmann} obtained the precise spectrum.
\begin{theorem}\label{Otto}
If $F$ is a finite field of characteristic $2$ and order $q$, then $\{W_{F,\card{F}-2}(a): a \in \Fu\}$ is the set $[1-2\sqrt{q},1+2\sqrt{q}]\cap (4 \Z)$.
\end{theorem}
And Katz and Livn\'e \cite{Katz-Livne} obtained the analogous result in characteristic $3$.
\begin{theorem}
If $F$ is a finite field of characteristic $3$ and order $q$, then $\{W_{F,\card{F}-2}(a): a \in \Fu\}$ is the set $[1-2\sqrt{q},1+2\sqrt{q}]\cap (3 \Z)$.
\end{theorem}
If the characteristic of $F$ is greater than $3$, Kononen, Rinta-aho, and V\"a\"an\"anen \cite{Kononen-Rinta-aho-Vaananen} showed a markedly different behavior.
\begin{theorem}
If $F$ is a finite field of characteristic $p > 3$, then there is no $a \in \Fu$ such that $W_{F,\card{F}-2}(a)=0$.
\end{theorem}
In characteristic $2$, a stronger version of the lower bound in Lemma \ref{Jacob} was deduced by Pursley and Sarwate \cite[eqs.~(4.6)--(4.8)]{Sarwate-Pursley} by applying a bound of of Sidel{\cprime}nikov \cite{Sidelnikov} to the family of Gold sequences (see Section \ref{Martha}) whose correlation spectra derive from $W_{F,d}$.
\begin{theorem}\label{Ignatius}
Let $F$ be a finite field of characteristic $2$ and $d$ an invertible exponent over $F$.
Then there is some $a \in \Fu$ such that $|W_{F,d}(a)-1| > \sqrt{2(\card{F}-2)}$.
If $[F:\F_2]$ is odd and greater than one, this means that there is some $a \in \Fu$ such that $|W_{F,d}(a)-1| \geq 1+\sqrt{2 \card{F}}$.
\end{theorem}
The $-1$ terms that appear in these bounds are a consequence of the fact that the correlation values are obtained from Weil sum values by subtracting $1$.
From these results about correlation, the following bound on the Weil sum itself immediately follows.
\begin{corollary}
Let $F$ be a finite field of characteristic $2$ with $[F:\F_2]$ odd, and let $d$ be an invertible exponent over $F$.
Then there is some $a \in \Fu$ such that $|W_{F,d}(a)| \geq \sqrt{2 \card{F}}$.
\end{corollary}
For every finite field $F$ of characteristic $2$ with $[F:\F_2]$ odd, there is some $d$ such that $W_{F,d}(a) \in \{0,\pm\sqrt{2 \card{F}}\}$ for all $a \in F$ (see the first, third, seventh, and ninth entries of Table \ref{Nancy} in Section \ref{Gerald} for examples), so the above corollary gives a sharp bound in this case.
When $[F:\F_2]$ is even, then Sarwate and Pursley \cite[\S3]{Sarwate-Pursley} suspect that the bound of Theorem \ref{Ignatius} is weak, and conjecture that it is not possible to keep all the values of $|W_{F,d}(a)-1|$ smaller than what one observes when $d=\card{F}-2$ in Theorem \ref{Otto}.
\begin{conjecture}[Sarwate-Pursley, 1980]
Let $F$ be of characteristic $2$ with $[F:\F_2]$ even, and let $d$ be an invertible exponent over $F$.
Then there is some $a \in \Fu$ such that $|W_{F,d}(a)-1| \geq -1 + 2\sqrt{\card{F}}$, or equivalently, there is some $a \in \Fu$ such that $|W_{F,d}(a)| \geq 2\sqrt{\card{F}}$.
\end{conjecture}
The equivalent version, not explicitly stated by Sarwate and Pursley, arises because every value $W_{F,d}(a)$ is a multiple of $4$ when $F$ is a non-prime field of characteristic $2$ (see Corollary \ref{Eugene}).
Feng, Leung, and Xiang \cite{Feng-Leung-Xiang} proved the following interesting result.
\begin{theorem}\label{Helen}
If $F$ is of characteristic $2$ and order $2^{2 m}$, then there is some $a \in \Fu$ such that $W_{F,d}(a) > 2^m+2^{\floor{m/2}}$.
\end{theorem}
In terms of the magnitudes $|W_{F,d}(a)|$, this is stronger than Theorem \ref{Ignatius} when $\card{F} \leq 64$.
But interestingly, the bound on $W_{F,d}(a)$ in Theorem \ref{Helen} does not involve an absolute value, so it tells us something specific about positive values of $W_{F,d}$.

\section{Non-Archimedean Bounds}\label{Felix}

The values $W_{F,d}(a)$ of our Weil sums lie in the cyclotomic field $\Q(\zeta_p)$, where $p$ is the characteristic of $F$ and $\zeta_p$ is a primitive $p$th root of unity.
The previous section dealt with Archimedean bounds, that is, bounds on the magnitude of our Weil sums with respect to the absolute value from the larger field $\C$.
We can also view our Weil sums $p$-adically, and ask for non-Archimedean bounds, that is, bounds on $p$-divisibility.
To this end, we recall that the {\it $p$-adic valuation} of a nonzero $a \in \Z$ is the unique $k$ such that $p^k \mid a$ and $p^{k+1}\nmid a$, where one defines $\valp(0)=\infty$, and then extends $\valp$ to $\Q$ by declaring that $\valp(a/b)=\valp(a)-\valp(b)$ for every $a,b \in \Z$ with $b\not=0$.
One should recall from algebraic number theory that $\Z[\zeta_p]$ is the ring of algebraic integers of the field $\Q(\zeta_p)$.
In $\Z[\zeta_p]$, the ideal $(p)$ is the $(p-1)$th power of the prime ideal $(1-\zeta_p)$.
This allows us to extend the usual $p$-adic valuation $\valp$ from $\Q$ to $\Q(\zeta_p)$, so that $v_p(1-\zeta_p)=1/(p-1)$.
We are interested in the $p$-adic valuations of Weil sum values $W_{F,d}(a)$.
\begin{definition}[Valuation of Weil Spectrum $V_{F,d}$]
Let $F$ be a finite field and let $d$ be invertible over $F$.
Then the {\it valuation of the Weil spectrum of $F$ and $d$}, denoted $V_{F,d}$, is
\[
V_{F,d} = \min_{a \in \Fu} \valp(W_{F,d}(a)).
\]
\end{definition}

To investigate $V_{F,d}$, it is useful to consider the Fourier expansion of the Weil sum in terms of Gauss sums.
If $F$ is a finite field of order $q$, we let $\zeta_{q-1}$ be a primitive $(q-1)$th root of unity over $\Q$, and then a {\it multiplicative character} is a group homomorphism $\chi\colon \Fu \to \Q(\zeta_{q-1})^*$.
We let $\mchars$ denote the group of multiplicative characters of $\Fu$, and $\id$ will always denote the trivial multiplicative character.
For any multiplicative character $\chi$ of a field $F$, the Gauss sum $G(\chi)$ is defined to be
\[
G(\chi)=\sum_{a \in \Fu} \psi(a)\chi(a),
\]
where $\psi$ is the canonical additive character of $F$.
Thus $G(\chi) \in \Q(\zeta_p,\zeta_{q-1})$.
One can expand the Weil sum in terms of pairwise products of these Gauss sums, and vice-versa.
\begin{lemma}
Let $F$ be a finite field of order $q$, let $d$ be an invertible exponent over $F$, and use $\id$ to denote the trivial multiplicative character of $F$.
Then 
\[
W_{F,d}(a) = \frac{q}{q-1} + \frac{1}{q-1} \sum_{\chi\not=\id} G(\chi) G(\conj{\chi}^d) \chi^d(-a),
\]
and
\[
\sum_{a \in \Fu} W_{F,d}(a) \conj{\chi^d(-a)} = \begin{cases}
q & \text{if $\chi=\id$}, \\
G(\chi) G(\conj{\chi}^d) & \text{otherwise.}
\end{cases}
\]
\end{lemma}
For a proof, see \cite[\S4]{Aubry-Katz-Langevin-15}, where the differences between the expressions there and here are due to a sign difference in the definition of $W_{F,d}(a)$.\footnote{And for the same reason, the formulae in \cite[Lemma 2.5]{Katz-Langevin-Lee-Sapozhnikov} should have the same sign changes, but one can see that this does not affect the proof of Corollary 2.6 of that paper, which is the only place this result is used.}
A consequence of our expansions is that we can understand the valuation of the Weil spectrum in terms of the valuations of pairwise products of Gauss sums.
\begin{corollary}\label{Clarence}
Let $F$ be a finite field of characteristic $p$, and let $d$ be an invertible exponent over $F$.  If $\card{F}=2$, then $d$ is degenerate over $F$ and $V_{F,d}=1$.  If $\card{F} > 2$, then
\[
V_{F,d}=\mins{\chi\in\mchars \\ \chi\not=\id} \valp(G(\chi) G(\conj{\chi}^d).
\]
\end{corollary}
See \cite[Corollary 2.6]{Katz-Langevin-Lee-Sapozhnikov} for a proof.\footnote{Also see \cite[Lemma 4.1]{Aubry-Katz-Langevin-15}, which neglects the fact that the trivial case when $\card{F}=2$ must be handled separately; the rest of that paper is unaffected by this oversight because the lemma is used only in Corollary 4.2 of that paper, which remains true because $V_{L,d} \leq [L:\F_2]$ can be deduced immediately from the fact that the first power moment for $W_{L,d}$ (given in Lemma 2.1(i) of that paper) has $2$-adic valuation $[L:\F_2]$.}
A further corollary of this result can be obtained using the Davenport-Hasse relation \cite[Corollary 4.2]{Aubry-Katz-Langevin-15}.
\begin{corollary}
Let $K$ be a finite field, $F$ a finite extension of $K$, and $d$ an invertible exponent over $F$.  Then $V_{F,d} \leq [F:K] V_{K,d}$.
\end{corollary}

Corollary \ref{Clarence} tells us how to obtain the valuation of the Weil spectrum from those of pairwise products of Gauss sums.
Now Stickelberger's Theorem tells us the $p$-adic valuation of these Gauss sums, which in turn gives the exact value of $V_{F,d}$ in terms of a combinatorial problem concerning quantities known as $p$-ary weights.
If $p$ is a prime and $n$ is an integer, we define the {\it $p$-ary weight function on $\Z/(p^n-1)\Z$}, to be the function $\wt \colon \Z/(p^n-1)\Z \to \Z$ with
\[
\wt(a_0 p^0 +a_1 p^1 + \cdots + a_{n-1} p^{n-1}) = a_0 + a_1 + \cdots + a_{n-1}
\]
where the powers of $p$ in this expression are elements of $\Z/(p^n-1)\Z$, while the coefficients $a_0,a_1,\ldots,a_{n-1}$ are elements of the subset $\{0,1,\ldots,p-1\}$ of $\Z$, with at least one $a_j < p-1$.
Since each element of $\Z/(p^n-1)\Z$ has a unique expression of this form, this is a well-defined function.
Now we may state what Stickelberger's Theorem tells us about the value of $V_{F,d}$.
\begin{lemma}
Let $F$ be a finite field of characteristic $p$ and order $q=p^n > 2$, let $d$ be an invertible exponent over $F$, and let $\wt$ be the $p$-ary weight function for $\Z/(p^n-1)\Z$.  Define
\[
m = \mins{j \in \Z/(q-1)\Z \\ j\not=0} \big(\wt(j)+\wt(-d j)\big),
\]
or equivalently
\[
m = (p-1) n + \mins{j \in \Z/(q-1)\Z \\ j\not=0} \big(\wt(d j)-\wt(j)\big).
\]
Then $V_{F,d} = m/(p-1)$.
\end{lemma}
When $d\equiv 1 \pmod{p-1}$ (which is invariably true in characteristic $2$), one can obtain an equivalent result from McEliece's Theorem \cite{McEliece-71,McEliece-72} on $p$-divisibility of weights in cyclic codes; earlier uses of this principle in characteristic $2$  often proceeded in this way \cite{McGuire-Calderbank,Calderbank-McGuire-Poonen-Rubinstein,Canteaut-Charpin-Dobbertin-99,Canteaut-Charpin-Dobbertin-00-Binary,Canteaut-Charpin-Dobbertin-00-Weight,Hollmann-Xiang,Hou,Dobbertin-Helleseth-Kumar-Martinsen}.  (But also see \cite{Langevin-Veron,Leander-Langevin} for approaches proceeding directly from Stickelberger's theorem in characteristic $2$.)
The more general result in the lemma here is stated and proved as \cite[Lemma 2.9]{Katz-Langevin-Lee-Sapozhnikov}.
The same result also appears as \cite[Proposition 4.3]{Katz-Langevin-15}, with the $\card{F} >2$ condition missing, which does not affect the other results of that paper.
One consequence of this result is that all values of our Weil sum have nontrivial $p$-divisibility, and we have some sense of the range of possible values for $V_{F,d}$.
\begin{corollary}\label{Eugene}
Let $F$ be a finite field of characteristic $p$ and let $d$ be an invertible exponent over $F$.
Then $V_{F,d}=[F:\Fp]$ if and only if $d$ is degenerate over $F$.
If $d$ is nondegenerate over $F$, then $2/(p-1) \leq V_{F,d} < [F:\Fp]$ with $V_{F,d}=2/(p-1)$ if and only if $d$ is an exponent equivalent to $\card{F}-2$ over $F$.
\end{corollary}
This was proved as \cite[Corollary 2.10]{Katz-Langevin-Lee-Sapozhnikov}.
Some stronger upper bounds for $V_{F,d}$ were also proved in \cite[Theorem 1.1, Remark 1.3]{Katz-Langevin-Lee-Sapozhnikov}.
\begin{theorem}\label{Julius}
Let $F$ be a finite field of characteristic $p$, and let $d$ be an invertible exponent over $F$.
\begin{enumerate}[(i)]
\item\label{Hans} If $d$ is degenerate over $F$, then $V_{F,d}=[F:\Fp]$.
\item\label{Edith} If $d$ is nondegenerate over $F$, but degenerate over $\Fp$, then we have the following:
\begin{enumerate}[(a)]
\item If\label{Ferdinand} $[F:\Fp]$ is a power of $2$, then $V_{F,d} \leq \frac{1}{2}[F:\Fp]$.
\item If\label{George} $[F:\Fp]$ is not a power of $2$, then $V_{F,d} \leq \frac{2}{3}[F:\Fp]$.
\end{enumerate}
\item\label{Henry} If $d$ is nondegenerate over $\Fp$ (which implies $p \geq 5$), then we have the following:
\begin{enumerate}[(a)]
\item\label{Irene} If $p\equiv 1\pmod{4}$ and $[F:\Fp]$ is odd, then $V_{F,d} \leq \frac{1}{2}[F:\Fp]$; but
\item\label{James} If $p\equiv 3\pmod{4}$ or $[F:\Fp]$ is even, then $V_{F,d} \leq \frac{1}{p-1} \ceil{\frac{p-1}{3}}[F:\Fp]$.
\end{enumerate}
\end{enumerate}
\end{theorem}
In many cases, the upper bounds of Theorem \ref{Julius} are the best possible, as discussed in \cite[Remark 1.3]{Katz-Langevin-Lee-Sapozhnikov}.
If we are not concerned about the degeneracy of the exponent over the prime field, we get the following simplified version of Theorem \ref{Julius}.
\begin{corollary}\label{Katherine}
Let $F$ be a finite field of characteristic $p$ and let $d$ be a nondegenerate invertible exponent over $F$.
Then $V_{F,d} \leq \frac{2}{3} [F:\Fp]$.
If $[F:\Fp]$ is a power of $2$, then $V_{F,d} \leq \frac{1}{2} [F:\Fp]$.
\end{corollary}
The following result \cite[Lemmata 4.1 and 4.2]{Katz-Langevin-Lee-Sapozhnikov} shows that the above bound is best possible when $[F:\Fp]$ is a power of $2$ greater than $1$.
\begin{theorem}
Let $F$ be a finite field of characteristic $p$ with $[F:\Fp]=2^s$ for some $s > 0$, and suppose that $F\not=\F_4$.
Then there is a nondegenerate invertible exponent $d$ over $F$ such that $V_{F,d}=\frac{1}{2}[F:\Fp]$.
\end{theorem}
Recall from Lemma \ref{Hubert} that the condition $F\not=\F_4$ is necessary for $F$ to have a nondegenerate exponent.
The general bound in Corollary \ref{Katherine} is best possible when $[F:\Fp]$ is divisible by $3$, as seen in the following result \cite[Lemma 3.2]{Katz-Langevin-Lee-Sapozhnikov}.
\begin{theorem}
Suppose $F$ is a finite field of characteristic $p$ with $[F:\Fp]$ not a power of $2$, and let $\ell$ be the smallest odd prime divisor of $[F:\Fp]$.
Then there is a nondegenerate invertible exponent $d$ over $F$ such that $V_{F,d}=\frac{\ell+1}{2\ell} [F:\Fp]$.
\end{theorem}
Katz, Langevin, Lee, and Sapozhnikov \cite[Conjecture 6.1]{Katz-Langevin-Lee-Sapozhnikov} conjecture that the $V_{F,d}$ values seen in this proposition are the highest possible.
\begin{conjecture}[Katz-Langevin-Lee-Sapozhnikov, 2017]\label{Lawrence}
If $F$ is a finite field of characteristic $p$ with $[F:\Fp]$ not a power of $2$, if $\ell$ is the smallest odd prime divisor of $[F:\Fp]$, and if $d$ is a nondegenerate invertible exponent over $F$, then $V_{F,d}\leq \frac{\ell+1}{2\ell} [F:\Fp]$.
\end{conjecture}
Conjecture \ref{Lawrence} proposes an upper bound that is often stronger than the $\frac{2}{3}[F:\Fp]$ bound of Corollary \ref{Katherine}, viz.~when $[F:\Fp]$ is greater than $10$ and is neither a power of $2$ nor a multiple of $3$.
Computer checks \cite[\S6]{Katz-Langevin-Lee-Sapozhnikov} have shown that Conjecture \ref{Lawrence} is true whenever $\card{F} < 10^{13}$.

\section{Number of Values}\label{Gerald}

In this section we are concerned with the cardinality of the Weil spectrum $\{W_{F,d}(a): a \in \Fu\}$.
\begin{definition}[Exactly, at Least, and at Most $k$-Valued $W_{F,d}$]
Let $F$ be a finite field and let $d$ be invertible over $F$.
We say that $W_{F,d}$ (or the Weil spectrum of $F$ and $d$) is {\it $k$-valued} (or sometimes {\it exactly $k$-valued} for emphasis) to mean that $\card{\{W_{F,d}(a): a \in \Fu\}}=k$.
And we say that $W_{F,d}$ is {\it at least $k$-valued} to mean that $\card{\{W_{F,d}(a): a \in \Fu\}} \geq k$, or {\it at most $k$-valued} to mean that $\card{\{W_{F,d}(a): a \in \Fu\}}\leq k$.
\end{definition}
If $d$ is degenerate over $F$, then Corollary \ref{Deidre} shows that $W_{F,d}$ is at most $2$-valued: to be precise, $W_{F,d}$ is $1$-valued when $F=\F_2$ and $2$-valued for any other field.
A classic result of Helleseth \cite[Theorem 4.1]{Helleseth-76} shows that $W_{F,d}$ is at least three-valued if $d$ is nondegenerate over $F$.
\begin{theorem}\label{Sally}
Let $F$ be a finite field and $d$ an invertible exponent over $F$.
Then $W_{F,d}$ is at least three-valued if and only if $d$ is nondegenerate.
\end{theorem}
There has been a great deal of interest in determining when $W_{F,d}$ is exactly three-valued.
Table \ref{Nancy} displays known pairs $(F,d)$ such that $W_{F,d}$ is three-valued; any known pair $(F,d)$ with a three-valued $W_{F,d}$ is either listed on the table or else corresponds to some entry $(F,d')$ on the table with $d'$ equivalent to $d$ over $F$.
Note also that we use the $2$-adic valuation $\valtwo$ in describing some of the required conditions on the table.
\begin{table}
\caption{Three-valued Weil sums of binomials}\label{Nancy}
\begin{center}
{\setlength{\extrarowheight}{3pt}\begin{tabular}{cccc}
order $q$ of $F$ & exponent $d$ & Weil spectrum & reference \\
\hline
 & $d=2^i+1$ & &  \\\noalign{\vskip-0.2pt}
$q=2^n$ & $0 < i < n$ & $0,\pm\sqrt{2^{\gcd(n,i)} q}$ & \cite{Kasami-66,Kasami-Lin-Peterson,Gold-68} \\\noalign{\vskip-0.2pt}
& $\valtwo(i) \geq \valtwo(n)$ & & \\\noalign{\vskip-0.2pt}
\rowcolor{tablecolor}[5.2pt]
 & $d=(p^{2 i}+1)/2$ & & \\\noalign{\vskip-0.2pt}
\rowcolor{tablecolor}[5.2pt]
\multirow{-2}{*}{$q=p^n$} & $0 < i < n$ & $0,\pm\sqrt{p^{\gcd(n,i)} q}$ & \multirow{-2}{*}{\cite{Trachtenberg} ($n$ odd)} \\\noalign{\vskip-0.2pt}
\rowcolor{tablecolor}[5.2pt]
\multirow{-2}{*}{$p$ odd} & $\valtwo(i) \geq \valtwo(n)$ &  & \multirow{-2}{*}{\cite{Helleseth-71,Helleseth-76} ($n$ even)}\\\noalign{\vskip-0.2pt}
 & $d=2^{2 i}-2^i+1$ & &  \\\noalign{\vskip-0.2pt}
$q=2^n$ & $0 < i < n$ & $0,\pm\sqrt{2^{\gcd(n,i)} q}$ & \cite{Welch,Kasami-71} \\\noalign{\vskip-0.2pt}
& $\valtwo(i) \geq \valtwo(n)$ & & \\\noalign{\vskip-0.2pt}
\rowcolor{tablecolor}[5.2pt]
 & $d=p^{2 i}-p^i+1$ & & \\\noalign{\vskip-0.2pt}
\rowcolor{tablecolor}[5.2pt]
\multirow{-2}{*}{$q=p^n$} & $0 < i < n$ & $0,\pm\sqrt{p^{\gcd(n,i)} q}$ & \multirow{-2}{*}{\cite{Trachtenberg} ($n$ odd)} \\\noalign{\vskip-0.2pt}
\rowcolor{tablecolor}[5.2pt]
\multirow{-2}{*}{$p$ odd} & $\valtwo(i) \geq \valtwo(n)$ &  & \multirow{-2}{*}{\cite{Helleseth-71,Helleseth-76} ($n$ even)}\\\noalign{\vskip-0.2pt}
$q=2^n$ & & & \\\noalign{\vskip-0.2pt}
$n>2$ & $d=2^{n/2}+2^{(n+2)/4}+1$ & $0,\pm 2\sqrt{q}$ & \cite{Cusick-Dobbertin} \\\noalign{\vskip-0.2pt}
$\valtwo(n)=1$ & & & \\\noalign{\vskip-0.2pt}
\rowcolor{tablecolor}[5.2pt]
$q=2^n$ & & & \\\noalign{\vskip-0.2pt}
\rowcolor{tablecolor}[5.2pt]
$n>2$ & $d=2^{(n+2)/2}+3$ & $0,\pm 2\sqrt{q}$ & \cite{Cusick-Dobbertin} \\\noalign{\vskip-0.2pt}
\rowcolor{tablecolor}[5.2pt]
$\valtwo(n)=1$ & & & \\\noalign{\vskip-0.2pt}
$q=2^n$ & & & \\\noalign{\vskip-0.2pt}
$n>1$ & $d=2^{(n-1)/2}+3$ & $0,\pm \sqrt{2 q}$ & \cite{Canteaut-Charpin-Dobbertin-99,Canteaut-Charpin-Dobbertin-00-Binary,Hollmann-Xiang} \\\noalign{\vskip-0.2pt}
$n$ odd & & & \\\noalign{\vskip-0.2pt}
\rowcolor{tablecolor}[5.2pt]
$q=3^n$ & & & \\\noalign{\vskip-0.2pt}
\rowcolor{tablecolor}[5.2pt]
$n>1$ & $d=2\cdot 3^{(n-1)/2}+1$ & $0,\pm \sqrt{3 q}$ & \cite{Dobbertin-Helleseth-Kumar-Martinsen} \\\noalign{\vskip-0.2pt}
\rowcolor{tablecolor}[5.2pt]
$n$ odd & & & \\\noalign{\vskip-0.2pt}
$q=2^n$ & $d=2^{2 i}+2^i-1$ & & \\\noalign{\vskip-0.2pt}
$n>1$ & $0 < i < n$ & $0,\pm \sqrt{2 q}$ & \cite{Hollmann-Xiang,Hou} \\\noalign{\vskip-0.2pt}
$n$ odd & $i \equiv -1/4 \pmod{n}$ & & \\\noalign{\vskip-0.2pt}
\rowcolor{tablecolor}[5.2pt]
$q=3^n$ & $d=2 \cdot 3^i+1$ & & \\\noalign{\vskip-0.2pt}
\rowcolor{tablecolor}[5.2pt]
$n>1$ & $0 < i <  n$ & $0,\pm \sqrt{3 q}$ & \cite{Katz-Langevin-15} \\\noalign{\vskip-0.2pt}
\rowcolor{tablecolor}[5.2pt]
$n$ odd & $i \equiv -1/4 \pmod{n}$  &  &
\end{tabular}}
\end{center}
\end{table}
Some features in Table \ref{Nancy} become immediately apparent.
First of all, it should be noted that for every entry, the three values of $W_{F,d}(a)$ lie in $\Z$.
That this must be true whenever $W_{F,d}$ is three-valued was proved in \cite[Theorem 1.7]{Katz-12}, and a simpler proof can be found in \cite[\S9]{Aubry-Katz-Langevin-15}.
\begin{theorem}\label{Agnes}
If $F$ is a finite field, $d$ is an invertible exponent over $F$, and $W_{F,d}$ is three-valued, then $d \equiv 1 \pmod{p-1}$ and $W_{F,d}(a) \in \Z$ for all $a \in F$.
\end{theorem}
Another feature of note in Table \ref{Nancy} is that $0$ always appears as one of the values of $W_{F,d}(a)$ with $a \in \Fu$.
This was proved in \cite[Theorem 1.9]{Katz-12}
\begin{theorem}\label{Norbert}
If $F$ is a finite field, $d$ is an invertible exponent over $F$, and $W_{F,d}$ is three-valued, then $W_{F,d}(a)=0$ for some $a \in \Fu$.
\end{theorem}
Another notable feature of the values of $W_{F,d}(a)$ appearing on the table is that the two nonzero values are opposites.
Lemma \ref{Julie} tells us that the two nonzero values must have opposite signs, but does not require them to have the same magnitude.
We say that a three-valued $W_{F,d}$ is {\it symmetric} if the two nonzero values in $\{W_{F,d}(a): a \in \Fu\}$ are opposites of each other.
It is not known whether all three-valued $W_{F,d}$ are symmetric, but we may conjecture it.
\begin{conjecture}\label{Simon}
If $F$ is a finite field, $d$ is an invertible exponent over $F$, and $W_{F,d}$ is three-valued, then it is symmetric.
\end{conjecture}
A much earlier conjecture looks at which fields $F$ have at least one invertible exponent $d$ such that $W_{F,d}$ is three-valued.
The first two entries of Table \ref{Nancy} show that for each field $F$ of characteristic $p$ and order $q=p^n$, one can obtain an exponent $d$ such that $W_{F,d}$ is three-valued, provided that $n$ is not a power of $2$.
And we note that when $n$ is a power of $2$, none of the entries of the table provide an exponent $d$ that makes $W_{F,d}$ three-valued.
This led Helleseth to conjecture (see \cite{Helleseth-71}, \cite[Conjecture 5.2]{Helleseth-76}) that fields of order $p^n$ with $n$ a power of $2$ do not support three-valued sums.
\begin{conjecture}[Helleseth Three-Valued Conjecture, 1971]\label{Thomas}
If $F$ is a finite field of characteristic $p$ with $[F:\Fp]$ a power of $2$, then there does not exist any invertible exponent $d$ over $F$ such that $W_{F,d}$ is three-valued.
\end{conjecture}
This conjecture has been studied by many researchers \cite{Calderbank-McGuire-Poonen-Rubinstein,McGuire,Charpin,CakCak-Langevin,Feng,Katz-12,Aubry-Katz-Langevin-14,Aubry-Katz-Langevin-15,Katz-15,Katz-Langevin-16}, and was proved to be true in characteristics $2$ and $3$ (see \cite[Corollary 1.10]{Katz-12} and \cite[Theorem 1.7]{Katz-15}).
\begin{theorem}\label{Alexandra}
If $F$ is a finite field of characteristic $p\in\{2,3\}$ with $[F:\Fp]$ a power of $2$, then there does not exist any invertible exponent $d$ over $F$ such that $W_{F,d}$ is three-valued.
\end{theorem}

Conjecture \ref{Thomas} remains open for characteristics $p \geq 5$, but some partial results have been obtained.
The following theorem (announced in \cite[Th\'eor\`eme 3.1]{Aubry-Katz-Langevin-14} and proved in \cite[Theorem 1.4]{Aubry-Katz-Langevin-15}) exhibits an obstruction to obtaining symmetric three-valued Weil sums.
\begin{theorem}\label{Linda}
Let $F$ be a finite field and $K$ and $L$ subfields of $F$ with $[L:K]=2$, and suppose that $d$ is an invertible exponent over $F$ that is nondegenerate over $L$ but degenerate over $K$.  Then $W_{F,d}$ is not symmetric three-valued.
\end{theorem}
A corollary of this theorem is that if there is any counterexample to Conjecture \ref{Thomas}, then it cannot be symmetric.
\begin{corollary}
If $F$ is a finite field of characteristic $p$ with $[F:\Fp]$ a power of $2$, and $d$ is an invertible exponent over $F$, then $W_{F,d}$ is not symmetric three-valued.
\end{corollary}
Thus Conjecture \ref{Simon} implies Conjecture \ref{Thomas}.
The characteristic $2$ case of this corollary was proved in \cite[Theorem 3]{Calderbank-McGuire-Poonen-Rubinstein}, and the general case was proved in \cite[Corollary 1.5]{Aubry-Katz-Langevin-15}.
A stronger version of the characteristic $2$ case due to Feng \cite[Theorem 2]{Feng} (stating that zero cannot appear in the Weil spectrum) was a critical advance toward proving Conjecture \ref{Thomas} in characteristic $2$; once this was established, the result that zero must occur (Theorem \ref{Norbert}) finished the proof of Helleseth's Three-Valued conjecture in characteristic $2$.

A symmetric three-valued Weil sum has many constraints on its three values, which we summarize here (see \cite[Proposition 1.3]{Aubry-Katz-Langevin-15}).
\begin{lemma}\label{David}
If $F$ is a finite field of characteristic $p$ and order $q$, and $d$ is an invertible exponent over $F$ such that $W_{F,d}$ is symmetric three-valued, then the three values must be $0$ and $\pm p^k$ for some positive integer $k$ such that $\sqrt{q} < p^k < q$.
\end{lemma}
A symmetric three-valued Weil sum is said to be {\it preferred} if the magnitude of the nonzero values is as small as Lemma \ref{David} allows it to be.
Thus when $F$ is of characteristic $p$ and order $q$, a preferred three-valued Weil sum has values $0$ and $\pm\sqrt{p q}$ if $q$ is an odd power of $p$, and has values $0$ and $\pm p\sqrt{q}$ if $q$ is an even power of $p$.
We note that the last six entries in Table \ref{Nancy} always provide preferred three-valued $W_{F,d}$ and the first four entries provide preferred three-valued $W_{F,d}$ for fields $F$ when $[F:\Fp]$ is not divisible by $4$ nor equal to $1$ or $2$.
In fact, the first two entries of Table \ref{Nancy} show that if $F$ has $[F:\Fp]$ not divisible by $4$ nor equal to $1$ or $2$, then there is some invertible exponent $d$ over $F$ such that $W_{F,d}$ is preferred three-valued.
The absence of preferred three-valued Weil sums when $[F:\Fp]$ is divisible by $4$ is is a consequence of a more general fact (see \cite[Th\'eor\`eme 3.3]{Aubry-Katz-Langevin-14} and \cite[Theorem 1.7]{Aubry-Katz-Langevin-15}) concerning the $2$-adic valuation of $[F:\Fp]$.
\begin{theorem}
Let $F$ be a finite field of characteristic $p$ and order $q$ with $\valtwo([F:\Fp])=s$.
If $d$ is an invertible exponent over $F$ such that $W_{F,d}$ is symmetric three-valued, then the three values must be $0$ and $\pm p^k$ for some positive integer $k$ such that $\sqrt{p^{2^s} q} \leq p^k < q$.
\end{theorem}
One consequence \cite[Corollary 1.8]{Aubry-Katz-Langevin-15} of this is that $W_{F,d}$ can never be preferred three-valued if the degree of $F$ over its prime field is divisible by $4$.
\begin{corollary}
If $F$ is a finite field of characteristic $p$ with $4\mid [F:\Fp]$, and $d$ is an invertible exponent over $F$, then $W_{F,d}$ is not preferred three-valued.
\end{corollary}

We saw in Theorem \ref{Linda} that there are constraints on three-valued $W_{F,d}$ when the exponent $d$ becomes degenerate in some proper subfield of $F$.
This prompts the following definition.
\begin{definition}[Niho Exponent]
If $F$ is a finite field of characteristic $p$ with $[F:\Fp]$ even, then we say that an invertible exponent $d$ over $F$ is a {\it Niho exponent over $F$} to mean that $d$ is nondegenerate over $F$ but $d$ is degenerate in the unique subfield $K$ of $F$ with $[F:K]=2$.
\end{definition}
These exponents are named after Niho, who studied them in \cite{Niho}.
It turns out that Niho exponents never yield three-valued Weil sums.
\begin{theorem}\label{Barbara}
If $F$ is a finite field and $d$ is a Niho exponent over $F$, then $W_{F,d}$ is not three-valued.
\end{theorem}
The characteristic $p=2$ case was proved in \cite[Theorem 2]{Charpin}, and the general result was proved in \cite[Theorem 9]{Helleseth-Lahtonen-Rosendahl}.
For a different proof in odd characteristic, see \cite[\S 8]{Aubry-Katz-Langevin-15}.
One consequence of Theorems \ref{Sally}, \ref{Agnes}, and \ref{Barbara} is that Conjecture \ref{Thomas} must be true in prime fields (see \cite[Corollary 1.8]{Katz-12}) and their quadratic extensions.
\begin{corollary}
If $F$ is a finite field of characteristic $p$ and order $p$ or $p^2$, and if $d$ is an invertible exponent over $F$, then $W_{F,d}$ is not three-valued.
\end{corollary}
Otherwise Theorem \ref{Agnes} would force $d$ to be degenerate in the prime field $\Fp$, which would make $d$ either degenerate or Niho over $F$, according to which $W_{F,d}$ could not be three-valued by Theorems \ref{Sally} and \ref{Barbara}.

Many of the results on three-valued $W_{F,d}$ use results on $p$-divisibility (see Section \ref{Felix}).
The following result \cite[Theorem 1.5]{Katz-15} shows that being three-valued constrains the valuation of $W_{F,d}$.
\begin{theorem}\label{Victor}
Let $F$ be a field of characteristic $p$ and order $q=p^n$, and suppose that $d$ is an invertible exponent such that $W_{F,d}$ is three-valued with values $0$, $a$, and $b$.  Then one of the following holds:
\begin{enumerate}[(i)]
\item $\valp(a) > n/2$ and $\valp(b) > n/2$; or
\item\label{Justin} $\valp(a)=\valp(b)=n/2$, and $|a-b|$ is a power of $p$ with $|a-b| > \sqrt{q}$;
\end{enumerate}
and case \eqref{Justin} cannot occur if $p\in\{2,3\}$.
\end{theorem}
This result leads to a proof of Theorem \ref{Alexandra}, which verifies Conjecture \ref{Thomas} in characteristics $2$ and $3$: we cannot have a three-valued $W_{F,d}$ with $F$ of characteristic $p=2$ or $3$ and $[F:\Fp]$ a power of $2$, because if we could, then Theorem \ref{Julius}\eqref{Ferdinand} would make $V_{F,d} \leq [F:\Fp]/2$ while Theorem \ref{Victor} would make $V_{F,d} > [F:\Fp]/2$.
This result also constrains the possibilities for a counterexample to Conjecture \ref{Thomas} in characteristics $p \geq 5$ (see \cite[Remark 4.6]{Katz-15} and \cite[\S5]{Katz-Langevin-16}).
Of course the main open problem which encompasses all others in the realm of three-valued $W_{F,d}$ is whether Table \ref{Nancy} is in fact a complete listing of all three-valued $W_{F,d}$.
\begin{openproblem}
Is every pair $(F,d)$ of finite field and invertible exponent such that $W_{F,d}$ is three-valued accounted for (up to the equivalence of Definition \ref{Elizabeth}) on Table \ref{Nancy}?
\end{openproblem}
Clearly if this were answered affirmatively, then Conjectures \ref{Simon} and \ref{Thomas} would also be answered affirmatively.

We have mainly focused on three-valued $W_{F,d}$ since three is the minimum number of values possible when $d$ is nondegenerate over $F$.
There is also interest in $(F,d)$ pairs such that $W_{F,d}$ takes on few values, albeit not as few as three.
Various examples of four-valued Weil spectra arising from Niho exponents have been found, and a synthesis and generalization of these results for fields of characteristic $2$ can be found in \cite[Theorem 23]{Dobbertin-Felke-Helleseth-Rosendahl}.

\section{Presence of a Zero Value}\label{Nigel}

In this section, we concern ourselves with the question as to when the Weil spectrum $\{W_{F,d}(a): a \in \Fu\}$ contains $0$.
If we had allowed $a$ to be $0$ in the definition of Weil spectrum, then this would be always be trivially true.
As a matter of fact, it is often the case that $W_{F,d}(a)=0$ for some nonzero $a$.
This prompted Helleseth to make the following conjecture (see \cite{Helleseth-71}, \cite[Conjecture 5.1]{Helleseth-76}).
\begin{conjecture}[Helleseth Vanishing Conjecture, 1971]\label{Rebecca}
If $F$ is a finite field of characteristic $p$ with $\card{F}>2$, and $d$ is an invertible exponent over $F$ with $d \equiv 1 \pmod{p-1}$, then $W_{F,d}(a)=0$ for some $a \in \Fu$.
\end{conjecture}
The original conjecture did not contain the proviso that $\card{F}>2$, which is necessary to exclude a trivial exception.
This conjecture appears to be even more difficult to approach then Helleseth's Three-Valued conjecture (presented as Conjecture \ref{Thomas} in this paper).
Theorem \ref{Norbert} can be taken as a proof of Conjecture \ref{Rebecca} in the special case where $W_{F,d}$ is three-valued.
An equivalent formulation of Conjecture \ref{Rebecca} can be expressed in terms of the {\it Dedekind determinant},
\[
D_{F,d} = \prod_{a\in \Fu} W_{F,d}(a),
\]
which leads to the following equivalent statement of Conjecture \ref{Rebecca}.
\begin{conjecture}[Helleseth 1971 Vanishing Conjecture, restated]
If $F$ is a finite field of characteristic $p$ with $\card{F}>2$, and $d$ is an invertible exponent over $F$ with $d \equiv 1 \pmod{p-1}$, then $D_{F,d}=0$.
\end{conjecture}
Although it has not been shown that $D_{F,d}$ always vanishes, some partial progress has been made by Aubry and Langevin \cite[Theorem 2]{Aubry-Langevin}, who show vanishing modulo certain primes.
\begin{theorem}
If $F$ is a finite field of characteristic $p$ and order $q=p^n > 2$, and $d$ is an invertible exponent over $F$ with $d \equiv 1 \pmod{p-1}$, then $D_{F,d} \equiv 0 \pmod{3}$.  Furthermore, if $n$ is a power of a prime $\ell$ with $\ell\nmid p-2$, then $D_{F,d} \equiv 0 \pmod{\ell}$.
\end{theorem}
Conjecture \ref{Rebecca} can also be expressed in terms of point counting.
If our finite field $F$ has order $q$, then Lemma \ref{Ursula} can be employed with $k=q-1$ and $\{b_1,\ldots,b_{q-1}\}=\Fu$ to calculate
\[
D_{F,d}=\frac{q^2 N_{b_1,\ldots,b_{q-1}}-q^{q-1}}{(q-1)^2},
\]
where $N_{b_1,\ldots,b_{q-1}}$ is as defined to be the number of solutions $(x_1,\ldots,x_{q-1}) \in F^{q-1}$ to the system of equations
\begin{align}
\begin{split}\label{Andrew}
b_1 x_1 + \cdots + b_{q-1} x_{q-1} & = 0 \\
x_1^d+\cdots +x_{q-1}^d & =0.
\end{split}
\end{align}
This gives the following restated version of Conjecture \ref{Rebecca}, proposed in \cite[p.~38]{Helleseth-02}.
\begin{conjecture}[Helleseth's 1971 Vanishing Conjecture, restated again]
Let $F$ be a finite field of characteristic $p$ and order $q > 2$, and let $d$ be an invertible exponent over $F$ with $d \equiv 1 \pmod{p-1}$.
Let the elements of $\Fu$ be written $b_1,\ldots,b_{q-1}$.
Then the number of solutions $(x_1,\ldots,x_{q-1}) \in F^{q-1}$ to the system of equations \eqref{Andrew} is equal to $q^{q-3}$.
\end{conjecture}
Thus the conjecture amounts to a precise determination of the cardinality of a hyperplane section of a Fermat variety.

\section*{Acknowledgement}

The author thanks Yakov Sapozhnikov, Arne Winterhof, and an anonymous reviewer, whose reading of and comments on the manuscript greatly helped the author.

\bibliographystyle{myabbrv}
\bibliography{181119a}

\begin{thebibliography}{10}

\bibitem{Aubry-Katz-Langevin-14}
Y.~Aubry, D.~J. Katz, and P.~Langevin.
\newblock Cyclotomie des sommes de {W}eil binomiales.
\newblock {\em C. R. Math. Acad. Sci. Paris}, {\bf 352}(5):373--376 (2014).

\bibitem{Aubry-Katz-Langevin-15}
Y.~Aubry, D.~J. Katz, and P.~Langevin.
\newblock Cyclotomy of {W}eil sums of binomials.
\newblock {\em J. Number Theory}, {\bf 154}:160--178 (2015).

\bibitem{Aubry-Langevin}
Y.~Aubry and P.~Langevin.
\newblock On a conjecture of {H}elleseth.
\newblock In {\em Algebraic informatics}, volume 8080 of {\em Lecture Notes in
  Comput. Sci.}, pages 113--118. Springer, Heidelberg (2013).

\bibitem{CakCak-Langevin}
E.~{\c{C}}ak{\c{C}}ak and P.~Langevin.
\newblock Power permutations in dimension 32.
\newblock In {\em Sequences and their applications---{SETA} 2010}, volume 6338
  of {\em Lecture Notes in Comput. Sci.}, pages 181--187. Springer, Berlin
  (2010).

\bibitem{Calderbank-McGuire-Poonen-Rubinstein}
A.~R. Calderbank, G.~McGuire, B.~Poonen, and M.~Rubinstein.
\newblock On a conjecture of {H}elleseth regarding pairs of binary
  {$m$}-sequences.
\newblock {\em IEEE Trans. Inform. Theory}, {\bf 42}(3):988--990 (1996).

\bibitem{Canteaut-Charpin-Dobbertin-99}
A.~Canteaut, P.~Charpin, and H.~Dobbertin.
\newblock Couples de suites binaires de longueur maximale ayant une
  corr\'elation crois\'ee \`a trois valeurs: conjecture de {W}elch.
\newblock {\em C. R. Acad. Sci. Paris S\'er. I Math.}, {\bf 328}(2):173--178
  (1999).

\bibitem{Canteaut-Charpin-Dobbertin-00-Binary}
A.~Canteaut, P.~Charpin, and H.~Dobbertin.
\newblock Binary {$m$}-sequences with three-valued crosscorrelation: a proof of
  {W}elch's conjecture.
\newblock {\em IEEE Trans. Inform. Theory}, {\bf 46}(1):4--8 (2000).

\bibitem{Canteaut-Charpin-Dobbertin-00-Weight}
A.~Canteaut, P.~Charpin, and H.~Dobbertin.
\newblock Weight divisibility of cyclic codes, highly nonlinear functions on
  {${\bf F}_{2^m}$}, and crosscorrelation of maximum-length sequences.
\newblock {\em SIAM J. Discrete Math.}, {\bf 13}(1):105--138 (2000).

\bibitem{Carlitz-78}
L.~Carlitz.
\newblock A note on exponential sums.
\newblock {\em Math. Scand.}, {\bf 42}(1):39--48 (1978).

\bibitem{Carlitz-79}
L.~Carlitz.
\newblock Explicit evaluation of certain exponential sums.
\newblock {\em Math. Scand.}, {\bf 44}(1):5--16 (1979).

\bibitem{Carlitz-Uchiyama}
L.~Carlitz and S.~Uchiyama.
\newblock Bounds for exponential sums.
\newblock {\em Duke Math. J.}, {\bf 24}:37--41 (1957).

\bibitem{Carmichael}
R.~D. Carmichael.
\newblock On sequences of integers defined by linear recurrence relations.
\newblock {\em Quarterly Journal of Pure and Applied Mathematics}, {\bf
  48}:343--372 (1920).

\bibitem{Charpin}
P.~Charpin.
\newblock Cyclic codes with few weights and {N}iho exponents.
\newblock {\em J. Combin. Theory Ser. A}, {\bf 108}(2):247--259 (2004).

\bibitem{Cochrane-Pinner-03}
T.~Cochrane and C.~Pinner.
\newblock Stepanov's method applied to binomial exponential sums.
\newblock {\em Q. J. Math.}, {\bf 54}(3):243--255 (2003).

\bibitem{Cochrane-Pinner-11}
T.~Cochrane and C.~Pinner.
\newblock Explicit bounds on monomial and binomial exponential sums.
\newblock {\em Q. J. Math.}, {\bf 62}(2):323--349 (2011).

\bibitem{Coulter}
R.~S. Coulter.
\newblock Further evaluations of {W}eil sums.
\newblock {\em Acta Arith.}, {\bf 86}(3):217--226 (1998).

\bibitem{Cusick-Dobbertin}
T.~W. Cusick and H.~Dobbertin.
\newblock Some new three-valued crosscorrelation functions for binary
  {$m$}-sequences.
\newblock {\em IEEE Trans. Inform. Theory}, {\bf 42}(4):1238--1240 (1996).

\bibitem{Davenport-Heilbronn}
H.~Davenport and H.~Heilbronn.
\newblock On an {E}xponential {S}um.
\newblock {\em Proc. London Math. Soc. (2)}, {\bf 41}(6):449--453 (1936).

\bibitem{Dobbertin-99-Niho}
H.~Dobbertin.
\newblock Almost perfect nonlinear power functions on {${\rm GF}(2^n)$}: the
  {N}iho case.
\newblock {\em Inform. and Comput.}, {\bf 151}(1-2):57--72 (1999).

\bibitem{Dobbertin-99-Welch}
H.~Dobbertin.
\newblock Almost perfect nonlinear power functions on {${\rm GF}(2^n)$}: the
  {W}elch case.
\newblock {\em IEEE Trans. Inform. Theory}, {\bf 45}(4):1271--1275 (1999).

\bibitem{Dobbertin-Felke-Helleseth-Rosendahl}
H.~Dobbertin, P.~Felke, T.~Helleseth, and P.~Rosendahl.
\newblock Niho type cross-correlation functions via {D}ickson polynomials and
  {K}loosterman sums.
\newblock {\em IEEE Trans. Inform. Theory}, {\bf 52}(2):613--627 (2006).

\bibitem{Dobbertin-Helleseth-Kumar-Martinsen}
H.~Dobbertin, T.~Helleseth, P.~V. Kumar, and H.~Martinsen.
\newblock Ternary {$m$}-sequences with three-valued cross-correlation function:
  new decimations of {W}elch and {N}iho type.
\newblock {\em IEEE Trans. Inform. Theory}, {\bf 47}(4):1473--1481 (2001).

\bibitem{Engstrom}
H.~T. Engstrom.
\newblock On sequences defined by linear recurrence relations.
\newblock {\em Trans. Amer. Math. Soc.}, {\bf 33}(1):210--218 (1931).

\bibitem{Feng}
T.~Feng.
\newblock On cyclic codes of length {$2^{2^r}-1$} with two zeros whose dual
  codes have three weights.
\newblock {\em Des. Codes Cryptogr.}, {\bf 62}(3):253--258 (2012).

\bibitem{Feng-Leung-Xiang}
T.~Feng, K.~Leung, and Q.~Xiang.
\newblock Binary cyclic codes with two primitive nonzeros.
\newblock {\em Sci. China Math.}, {\bf 56}(7):1403--1412 (2013).

\bibitem{Games-m-sequences}
R.~A. Games.
\newblock The geometry of {$m$}-sequences: three-valued crosscorrelations and
  quadrics in finite projective geometry.
\newblock {\em SIAM J. Algebraic Discrete Methods}, {\bf 7}(1):43--52 (1986).

\bibitem{Games-quadrics}
R.~A. Games.
\newblock The geometry of quadrics and correlations of sequences.
\newblock {\em IEEE Trans. Inform. Theory}, {\bf 32}(3):423--426 (1986).

\bibitem{Gold-67}
R.~Gold.
\newblock Optimal binary sequences for spread spectrum multiplexing.
\newblock {\em IEEE Transactions on Information Theory}, {\bf 13}(4):619--621
  (1967).

\bibitem{Gold-68}
R.~Gold.
\newblock Maximal recursive sequences with 3-valued recursive cross-correlation
  functions.
\newblock {\em IEEE Transactions on Information Theory}, {\bf 14}(1):154--156
  (1968).

\bibitem{Golomb}
S.~W. Golomb.
\newblock {\em Shift register sequences}.
\newblock With portions co-authored by Lloyd R. Welch, Richard M. Goldstein,
  and Alfred W. Hales. Holden-Day, Inc., San Francisco,
  Calif.-Cambridge-Amsterdam (1967).

\bibitem{Golomb-Gong}
S.~W. Golomb and G.~Gong.
\newblock {\em Signal design for good correlation}.
\newblock Cambridge University Press, Cambridge (2005).

\bibitem{Hall}
M.~Hall.
\newblock An isomorphism between linear recurring sequences and algebraic
  rings.
\newblock {\em Trans. Amer. Math. Soc.}, {\bf 44}(2):196--218 (1938).

\bibitem{Helleseth-71}
T.~Helleseth.
\newblock Krysskorrelasjonsfunksjonen mellom maksimale sekvenser over {GF(q)}.
\newblock Master's thesis, Matematisk Institutt, Universitetet i Bergen (1971).

\bibitem{Helleseth-76}
T.~Helleseth.
\newblock Some results about the cross-correlation function between two maximal
  linear sequences.
\newblock {\em Discrete Math.}, {\bf 16}(3):209--232 (1976).

\bibitem{Helleseth-02}
T.~Helleseth.
\newblock On the crosscorrelation of {$m$}-sequences and related sequences with
  ideal autocorrelation.
\newblock In {\em Sequences and their applications ({B}ergen, 2001)}, Discrete
  Math. Theor. Comput. Sci. (Lond.), pages 34--45. Springer, London (2002).

\bibitem{Helleseth-Lahtonen-Rosendahl}
T.~Helleseth, J.~Lahtonen, and P.~Rosendahl.
\newblock On {N}iho type cross-correlation functions of {$m$}-sequences.
\newblock {\em Finite Fields Appl.}, {\bf 13}(2):305--317 (2007).

\bibitem{Hollmann-Xiang}
H.~D.~L. Hollmann and Q.~Xiang.
\newblock A proof of the {W}elch and {N}iho conjectures on cross-correlations
  of binary {$m$}-sequences.
\newblock {\em Finite Fields Appl.}, {\bf 7}(2):253--286 (2001).

\bibitem{Hou}
X.-D. Hou.
\newblock A note on the proof of {N}iho's conjecture.
\newblock {\em SIAM J. Discrete Math.}, {\bf 18}(2):313--319 (2004).

\bibitem{Karatsuba}
A.~A. Karatsuba.
\newblock Estimates of complete trigonometric sums.
\newblock {\em Mat. Zametki}, {\bf 1}(2):199--208 (1967).

\bibitem{Kasami-66}
T.~Kasami.
\newblock Weight distribution formula for some class of cyclic codes.
\newblock Technical report, Coordinated Science Laboratory, University of
  Illinois, Urbana (1966).

\bibitem{Kasami-71}
T.~Kasami.
\newblock The weight enumerators for several classes of subcodes of the {$2$}nd
  order binary {R}eed-{M}uller codes.
\newblock {\em Information and Control}, {\bf 18}:369--394 (1971).

\bibitem{Kasami-Lin-Peterson}
T.~Kasami, S.~Lin, and W.~W. Peterson.
\newblock Some results on cyclic codes which are invariant under the affine
  group and their applications.
\newblock {\em Information and Control}, {\bf 11}:475--496 (1967).

\bibitem{Katz-12}
D.~J. Katz.
\newblock Weil sums of binomials, three-level cross-correlation, and a
  conjecture of {H}elleseth.
\newblock {\em J. Combin. Theory Ser. A}, {\bf 119}(8):1644--1659 (2012).

\bibitem{Katz-15}
D.~J. Katz.
\newblock Divisibility of {W}eil sums of binomials.
\newblock {\em Proc. Amer. Math. Soc.}, {\bf 143}(11):4623--4632 (2015).

\bibitem{Katz-Langevin-15}
D.~J. Katz and P.~Langevin.
\newblock Proof of a conjectured three-valued family of {W}eil sums of
  binomials.
\newblock {\em Acta Arith.}, {\bf 169}(2):181--199 (2015).

\bibitem{Katz-Langevin-16}
D.~J. Katz and P.~Langevin.
\newblock New open problems related to old conjectures by {H}elleseth.
\newblock {\em Cryptogr. Commun.}, {\bf 8}(2):175--189 (2016).

\bibitem{Katz-Langevin-Lee-Sapozhnikov}
D.~J. Katz, P.~Langevin, S.~Lee, and Y.~Sapozhnikov.
\newblock The {$p$}-adic valuations of {W}eil sums of binomials.
\newblock {\em J. Number Theory}, {\bf 181}:1--26 (2017).

\bibitem{Katz-Livne}
N.~Katz and R.~Livn\'e.
\newblock Sommes de {K}loosterman et courbes elliptiques universelles en
  caract\'eristiques {$2$} et {$3$}.
\newblock {\em C. R. Acad. Sci. Paris S\'er. I Math.}, {\bf 309}(11):723--726
  (1989).

\bibitem{Kononen-Rinta-aho-Vaananen}
K.~P. Kononen, M.~J. Rinta-aho, and K.~O. V\"a\"an\"anen.
\newblock On integer values of {K}loosterman sums.
\newblock {\em IEEE Trans. Inform. Theory}, {\bf 56}(8):4011--4013 (2010).

\bibitem{Lachaud-Wolfmann}
G.~Lachaud and J.~Wolfmann.
\newblock Sommes de {K}loosterman, courbes elliptiques et codes cycliques en
  caract\'eristique {$2$}.
\newblock {\em C. R. Acad. Sci. Paris S\'er. I Math.}, {\bf 305}(20):881--883
  (1987).

\bibitem{Langevin-Veron}
P.~Langevin and P.~V\'eron.
\newblock On the non-linearity of power functions.
\newblock {\em Des. Codes Cryptogr.}, {\bf 37}(1):31--43 (2005).

\bibitem{Leander-Langevin}
G.~Leander and P.~Langevin.
\newblock On exponents with highly divisible {F}ourier coefficients and
  conjectures of {N}iho and {D}obbertin.
\newblock In {\em Algebraic geometry and its applications}, volume~5 of {\em
  Ser. Number Theory Appl.}, pages 410--418. World Sci. Publ., Hackensack, NJ
  (2008).

\bibitem{McEliece-71}
R.~J. McEliece.
\newblock On periodic sequences from {${\rm GF}(q)$}.
\newblock {\em J. Combinatorial Theory Ser. A}, {\bf 10}:80--91 (1971).

\bibitem{McEliece-72}
R.~J. McEliece.
\newblock Weight congruences for {$p$}-ary cyclic codes.
\newblock {\em Discrete Mathematics}, {\bf 3}(1):177 -- 192 (1972).

\bibitem{McGuire}
G.~McGuire.
\newblock On certain 3-weight cyclic codes having symmetric weights and a
  conjecture of {H}elleseth.
\newblock In {\em Sequences and their applications ({B}ergen, 2001)}, Discrete
  Math. Theor. Comput. Sci. (Lond.), pages 281--295. Springer, London (2002).

\bibitem{McGuire-Calderbank}
G.~McGuire and A.~R. Calderbank.
\newblock Proof of a conjecture of {S}arwate and {P}ursley regarding pairs of
  binary {$m$}-sequences.
\newblock {\em IEEE Trans. Inform. Theory}, {\bf 41}(4):1153--1155 (1995).

\bibitem{Mordell}
L.~{Mordell}.
\newblock {On a sum analogous to a Gauss's sum.}
\newblock {\em {Q. J. Math., Oxf. Ser.}}, {\bf 3}:161--167 (1932).

\bibitem{Niho}
Y.~Niho.
\newblock {\em Multi-valued cross-correlation function between two maximal
  linear recursive sequences}.
\newblock PhD thesis, University of Southern California, Los Angeles (1972).

\bibitem{Sarwate-Pursley}
D.~V. Sarwate and M.~B. Pursley.
\newblock Crosscorrelation properties of pseudorandom and related sequences.
\newblock {\em Proceedings of the IEEE}, {\bf 68}(5):593--619 (1980).
\newblock Correction in {\em Proceedings of the IEEE}, {\bf 68}(12):1554
  (1980).

\bibitem{Sidelnikov}
V.~M. Sidel{\cprime}nikov.
\newblock The mutual correlation of sequences.
\newblock {\em Dokl. Akad. Nauk SSSR}, {\bf 196}:531--534 (1971).
\newblock English translation in {\em Soviet Math. Dokl.}, {\bf 12}:197--201
  (1971).

\bibitem{Trachtenberg}
H.~M. Trachtenberg.
\newblock {\em On the cross-correlation functions of maximal linear sequences}.
\newblock PhD thesis, University of Southern California, Los Angeles (1970).

\bibitem{Ward-31}
M.~Ward.
\newblock The distribution of residues in a sequence satisfying a linear
  recursion relation.
\newblock {\em Trans. Amer. Math. Soc.}, {\bf 33}(1):166--190 (1931).

\bibitem{Ward-33}
M.~Ward.
\newblock The arithmetical theory of linear recurring series.
\newblock {\em Trans. Amer. Math. Soc.}, {\bf 35}(3):600--628 (1933).

\bibitem{Weil}
A.~Weil.
\newblock On some exponential sums.
\newblock {\em Proc. Nat. Acad. Sci. U. S. A.}, {\bf 34}:204--207 (1948).

\bibitem{Welch}
L.~R. Welch.
\newblock Trace mappings in finite fields and shift register cross-correlation
  properties.
\newblock Technical report, Dept. Electrical Engineering, University of
  Southern California, Los Angeles (1969).

\end{thebibliography}
\end{document}